\documentclass[11pt,reqno,a4paper]{amsart}
\usepackage[T1]{fontenc}
\usepackage{lmodern,microtype,amsmath,amssymb}
\usepackage[margin=28mm]{geometry}
\usepackage{booktabs,tikz,placeins,capt-of}
\usepackage[hidelinks]{hyperref}
\hypersetup{pdftitle={Torus actions and duality for exponential polynomial solutions of linear differential equations}}
\newtheorem{theorem}{Theorem}[section]
\newtheorem{lemma}[theorem]{Lemma}
\newtheorem{corollary}[theorem]{Corollary}
\newtheorem{proposition}[theorem]{Proposition}
\theoremstyle{remark}
\newtheorem*{remark}{Remark}
\numberwithin{equation}{section}
\newcommand{\C}{\mathbb C}
\newcommand{\Z}{\mathbb Z}
\newcommand{\R}{\mathbb R}
\newcommand{\Cs}{\C^*}
\DeclareMathOperator{\conv}{conv}
\title[Torus actions and duality]{Torus actions and duality for exponential polynomial\\
solutions of linear differential equations}
\author{Xing-Yu Li}
\address{\textnormal{Xing-Yu Li}\newline
\hspace*{\parindent}Department of Mathematics, The University of Manchester,
Oxford Road, Manchester M13~9PL, United Kingdom}
\email{\href{mailto:xingyu.li@manchester.ac.uk}{xingyu.li@manchester.ac.uk}}
\urladdr{\href{https://xyli9596.github.io/}{https://xyli9596.github.io/}}
\thanks{The author acknowledges support from the China Scholarship Council
(CSC No.~202508440485).}
\date{}
\begin{document}
\raggedbottom
\begin{abstract}
For \(f''+Af'+Bf=0\) with exponential polynomial coefficients, we
prove that any two exponents of a nonzero exponential polynomial
solution differ by a polynomial of degree at most one. A limit of
the torus action isolates the terms on an edge of the Newton polytope,
where a ramification argument yields the bound. Building on earlier work
of Heittokangas, Ishizaki, Tohge and Wen, we use this bound to show that,
in their setting, the common growth order of \(f\) and \(A\) is one and
\(B\) is a polynomial, answering one of their questions. We also give a
criterion for the least order of a monic linear differential equation
that has exponential polynomial coefficients and admits a given nonzero
exponential polynomial as a solution. For each \(n\ge3\), we construct
equations of order \(n\) in which a solution and
the unique dominant coefficient are dual and have arbitrarily large
common growth order, answering another question posed by the same authors.
\end{abstract}
\maketitle
\enlargethispage{6pt}
\vspace{-24pt}

\section{Introduction}\label{sec:introduction}
An \emph{exponential polynomial} is an entire function of the form
\begin{equation}\label{eq:sum}
 f(z)=\sum_{j=1}^m P_j(z)e^{Q_j(z)},\qquad P_j,Q_j\in\C[z].
\end{equation}
An exponent \(Q\) is \emph{normalized} if it vanishes at zero. We use
distinct normalized \(Q_j\) and nonzero \(P_j\).
If \(f\) is transcendental, its growth order is \(q=\max_j\deg Q_j\),
and the coefficients of \(z^q\) in the \(Q_j\) are its
\emph{highest-order frequencies}.

Two transcendental exponential polynomials of the same growth order
are called \emph{dual} if the nonzero highest-order frequencies of one
lie on a ray from the origin and those of the other lie on the opposite ray
\cite[Definition~4.12]{survey}.

Consider the linear differential equation
\begin{equation}\label{eq:general}
 f^{(n)}+\sum_{j=0}^{n-1}A_j f^{(j)}=0,
 \qquad A_j\text{ entire}.
\end{equation}
If \(j\) is the largest index for which \(A_j\) is transcendental,
Frei \cite[p.~207, (2.2)]{frei-growth} proved that \eqref{eq:general}
can have at most \(j\) linearly independent solutions of finite order.

Wen, Gundersen and Heittokangas \cite[Theorem~2]{wgh} proved that if
the \(A_j\) are exponential polynomials and a transcendental coefficient
\(A_\mu\), with \(1\le\mu<n\), has strictly larger growth order than
every other coefficient, then \(A_\mu\) is dual to every transcendental
exponential polynomial solution of \eqref{eq:general}.
For example, \(f=1+e^z\) solves Frei's equation
\(f''+e^{-z}f'-f=0\) \cite{frei}. The solution and the dominant
coefficient \(e^{-z}\) are dual of growth order one, with nonzero
frequencies \(1\) and \(-1\), respectively.

For the second-order equation
\begin{equation}\label{eq:ode}
 f''+Af'+Bf=0,
\end{equation}
Heittokangas et al.\ \cite[Theorem~3.1]{dual} considered transcendental
exponential polynomials \(f,A\) and an entire function \(B\) with
\(T(r,B)=o(T(r,A))\), where \(T\) is the Nevanlinna characteristic.
They showed that \(f\) and \(A\) are dual of some common growth order
\(q\), and that \(B\) is an exponential polynomial of order at most
\(q-1\). They asked whether the common order \(q\) must be one and
\(B\) must be a polynomial \cite[Problem~1]{dual}, and restated the
question in their survey
\cite[Problem~6]{survey}.

The affirmative answer follows from a restriction on the differences
\(Q_i-Q_j\) between the exponents of a solution. This restriction
holds for arbitrary exponential polynomial coefficients \(A\) and
\(B\), regardless of duality.

\begin{theorem}\label{thm:exponents}
Let \(A\) and \(B\) be exponential polynomials, and let a nonzero
solution of \eqref{eq:ode} have the form \eqref{eq:sum}.
Then \(Q_i-Q_j\) has degree at most one for every pair \(i,j\).
Moreover, \(f\) has \(O(r)\) zeros in \(|z|\le r\), counted with
multiplicity.
\end{theorem}

The degree bound gives \(Q_j(z)-Q_1(z)=a_jz\) for distinct constants
\(a_j\), since the \(Q_j\) are normalized. Thus
\(f=e^{Q_1}g\), where \(g(z)=\sum_{j=1}^m P_j(z)e^{a_jz}\) has the
same zeros as \(f\), with multiplicities. If \(m\ge2\), Dickson's
theorem \cite{dickson} gives \(n_f(r)=cr+O(1)\), where \(n_f(r)\)
counts the zeros in \(|z|\le r\) and \(c>0\) depends only on the
\(a_j\). The zeros then have exponent of convergence one. If
\(m=1\), all the zeros are those of the polynomial \(P_1\).

Both conclusions of Theorem~\ref{thm:exponents} can fail for
third-order equations. With \(H=1\),
\cite[Example~4.1]{dual} gives
\[
 f'''-2zf''+\bigl(e^{-z^2}-3\bigr)f'-2zf=0
\]
with solution \(f=1+e^{z^2}\). The coefficient of \(f'\) is the unique
transcendental coefficient and is dual to \(f\). The exponents of \(f\)
differ by \(z^2\), and \(f\) has \(2r^2/\pi+O(1)\) zeros in \(|z|\le r\).

Under the hypotheses of the question above, the bound on exponent
differences also controls the common growth order of \(f\) and \(A\).

\begin{corollary}\label{thm:hitw}
Suppose that \(f\) and \(A\) in \eqref{eq:ode} are transcendental
exponential polynomials and that \(B\) is entire with
\(T(r,B)=o(T(r,A))\). Then \(\rho(f)=\rho(A)=1\) and \(B\)
is a polynomial.
\end{corollary}

The third-order example above attains the common growth order
\(q=n-1=2\). For \eqref{eq:general}, Heittokangas et al.\
\cite[Problem~2]{dual} asked whether the common growth order \(q\)
of a solution and its dual dominant coefficient must satisfy
\(q\le n-1\). The following construction answers this question in
the negative.

\begin{proposition}\label{prop:higher}
For each pair of integers \(n\ge3\) and \(q\ge2\), there is an
equation \eqref{eq:general} with a solution \(f\) such that \(f\) and
the coefficient \(A_1\) of \(f'\) are dual exponential polynomials of
growth order \(q\). All other coefficients are polynomials.
\end{proposition}

Sections~\ref{sec:family} and~\ref{sec:branching} prove
Theorem~\ref{thm:exponents} by passing from a family of twisted
solutions to a limit supported on an edge of the Newton polytope.
Section~\ref{sec:hitw} proves Corollary~\ref{thm:hitw}.
Section~\ref{sec:prescribed} uses the same twists to give a criterion
for solving an equation of prescribed order and proves
Proposition~\ref{prop:higher}.

The Lean 4 supplement verifies the statements of Theorem~\ref{thm:exponents}
and Propositions~\ref{prop:higher} and~\ref{prop:twist-order}.
It also verifies the derivation of Corollary~\ref{thm:hitw} with the
structure theorem of Heittokangas et al.\ \cite[Theorem~4.16(a)]{survey}
as an explicit hypothesis.

\section{A family of equations with simple zeros}\label{sec:family}

Multiplying \(e^z\) by \(v\ne0\) requires multiplying \(e^{-z}\)
by \(v^{-1}\) to preserve \(e^ze^{-z}=1\). We use characters to
rescale the terms of \(f,A\) and \(B\) and obtain a family of solutions
whose zeros are all simple.

\subsection{The group of exponents and its characters}
Let \(\Gamma\) consist of all integer linear combinations of the
exponents of \(f,A\) and \(B\), normalized as in Section~\ref{sec:introduction}.
It is an additive group of polynomials. A \(\Z\)-basis
\(\gamma_1,\ldots,\gamma_r\) gives each element of \(\Gamma\) a unique
expression \(\sum_{i=1}^r m_i\gamma_i\) with \(m_i\in\Z\). An additive
group with such a basis is called
\emph{free abelian}. In these coordinates \(\Gamma\) is the integer
lattice \(\Z^r\) inside \(\R^r\).

For example, suppose the exponents collected from \(f,A\) and \(B\) are
\[
 z^2,\qquad \frac{z}{3},\qquad z^2+iz,\qquad z^2+\sqrt{2}\,z,
 \qquad \frac{z^2+z}{2}.
\]
In the \(\mathbb Q\)-basis \(z^2,z,iz,\sqrt{2}\,z\) of their rational
span, the second and fifth vectors have coordinates \((0,1/3,0,0)\)
and \((1/2,1/2,0,0)\). Multiplying all coordinates by \(6\) sends
these vectors to \((0,2,0,0)\) and \((3,3,0,0)\) and gives all five
generators integer coordinates. A \(\Z\)-basis of \(\Gamma\) is
\((z^2+z)/2,z/3,iz,\sqrt{2}\,z\).

In general, clearing a common denominator in the rational coordinates
of the generators embeds \(\Gamma\) in \(\Z^d\), where \(d\) is the
dimension of their rational span. Every subgroup of \(\Z^d\) is free
abelian \cite[Proposition~A.24]{hatcher}, so \(\Gamma\) has a \(\Z\)-basis.

\noindent\textbf{Character.}\enspace
A \emph{character} of \(\Gamma\) is a group homomorphism
\(\chi\colon\Gamma\to\Cs\) from the additive group \(\Gamma\) to the
multiplicative group \(\Cs=\C\setminus\{0\}\), so that
\(\chi(\gamma+\delta)=\chi(\gamma)\chi(\delta)\). The identity gives
\(\chi(0)=1\) and \(\chi(-\gamma)=\chi(\gamma)^{-1}\).
At each basis element \(\gamma_i\), we can prescribe the value
\(w_i=\chi(\gamma_i)\) as any nonzero complex number. The value at
\(\gamma=\sum_i m_i\gamma_i\) is then \(\prod_i w_i^{m_i}\).
For example, in \(\Z z\oplus\Z z^2\) a character takes the form
\(\chi(mz+nz^2)=v^m w^n\) for \(v,w\in\Cs\).
These integer powers also appear in the usual description of
characters and one-parameter subgroups \cite[Lecture~I, Sect.~2]{cox}.

\subsection{Twisting the equation}\label{subsec:twisting}
\leavevmode\par
\noindent\textbf{Torus action.}\enspace
The group \((\Cs)^r\) with componentwise multiplication is a
\emph{complex algebraic torus}. A \emph{torus action} on a set \(V\) is
a map \((\Cs)^r\times V\to V\), written \((\mathbf w,v)\mapsto\mathbf w\cdot v\),
satisfying
\[
 \mathbf 1\cdot v=v,\qquad
 \mathbf u\cdot(\mathbf w\cdot v)=(\mathbf u\mathbf w)\cdot v
\]
for all \(v\in V\) and \(\mathbf u,\mathbf w\in(\Cs)^r\), where
\(\mathbf 1=(1,\ldots,1)\). A Laurent polynomial in \(w_1,\ldots,w_r\)
is a finite linear combination of Laurent monomials
\(w_1^{m_1}\cdots w_r^{m_r}\) with \(m_i\in\Z\).
On a finite-dimensional complex vector space \(V\), the action is
\emph{algebraic} when its components in linear coordinates are
polynomials in the coordinates of \(v\), with Laurent polynomial
coefficients in \(w_1,\ldots,w_r\) \cite[pp.~10--11 and~199]{cls}.

\noindent\textbf{Character twist.}\enspace
We write \(F=\sum_\gamma P_\gamma e^\gamma\) for any exponential
polynomial with exponents in \(\Gamma\) and arbitrary polynomial
multipliers \(P_\gamma\), and will apply the construction below to
\(f,A\) and \(B\) in turn. Multiplying each term \(P_\gamma e^\gamma\)
by the constant \(\chi(\gamma)\) assigned to its exponent gives the
\emph{twist}
\begin{equation}\label{eq:twist-definition}
 F_\chi=\sum_\gamma P_\gamma\chi(\gamma)e^\gamma.
\end{equation}
The exponentials \(e^\gamma\) with distinct normalized exponents are
linearly independent over \(\C(z)\) \cite[Sect.~3, after (3.1)]{survey}.
Thus the representation of \(F\) is unique and its twist is well defined.
For \(\mathbf w=(w_1,\ldots,w_r)\in(\Cs)^r\), let \(\chi_{\mathbf w}\)
be the unique character with \(\chi_{\mathbf w}(\gamma_i)=w_i\) for each
\(i\) and define \(\mathbf w\cdot F=F_{\chi_{\mathbf w}}\) by taking
\(\chi=\chi_{\mathbf w}\) in \eqref{eq:twist-definition}.
Coordinatewise multiplication gives \(\chi_{\mathbf 1}(\gamma)=1\) and
\(\chi_{\mathbf u\mathbf w}(\gamma)
=\chi_{\mathbf u}(\gamma)\chi_{\mathbf w}(\gamma)\).
Thus \eqref{eq:twist-definition} gives \(\mathbf 1\cdot F=F\) and
\(\mathbf u\cdot(\mathbf w\cdot F)=(\mathbf u\mathbf w)\cdot F\).
This verifies the compatibility condition for a group action.

For a fixed \(F\), its nonzero terms form a basis of the space
\(V_F=\operatorname{span}_{\C}\{P_\gamma e^\gamma:P_\gamma\not\equiv0\}\).
For \(H=\sum_\gamma a_\gamma P_\gamma e^\gamma\in V_F\), we have
\[
 \mathbf w\cdot H
 =\sum_\gamma a_\gamma\left(\prod_{i=1}^r w_i^{m_i(\gamma)}\right)
 P_\gamma e^\gamma\in V_F,
 \qquad \gamma=\sum_{i=1}^r m_i(\gamma)\gamma_i.
\]
Each coordinate \(a_\gamma\) is multiplied by a Laurent monomial in
\(\mathbf w\), so the action on \(V_F\) is algebraic.

When two terms with exponents \(\gamma\) and \(\delta\) are multiplied,
their exponents add. The defining property of \(\chi\) therefore gives
\((FG)_\chi=F_\chi G_\chi\). Differentiating a term gives
\((P_\gamma e^\gamma)'=(P_\gamma'+\gamma'P_\gamma)e^\gamma\), so
\((F')_\chi=(F_\chi)'\). Applying \eqref{eq:twist-definition} separately
to \(f,A\) and \(B\) gives \(f_\chi,A_\chi,B_\chi\), respectively.
By linearity and the product and differentiation rules, \eqref{eq:ode}
then yields
\begin{equation}\label{eq:twist}
 (f''+Af'+Bf)_\chi=f_\chi''+A_\chi f_\chi'+B_\chi f_\chi=0.
\end{equation}
For example, Frei's equation \(f''+e^{-z}f'-f=0\) has the solution
\(f=1+e^z\). Here \(\Gamma=\Z z\), and \(w\in\Cs\) determines the
character \(\chi_w\) with \(\chi_w(z)=w\). Its values at \(0,z,-z\)
are \(1,w,w^{-1}\), respectively. Applying \eqref{eq:twist-definition}
to \(f,A\) and \(B\) with this character gives \(f_w=1+we^z\),
\(A_w=w^{-1}e^{-z}\), and \(B_w=-1\). The twisted equation is
\[
 f_w''+w^{-1}e^{-z}f_w'-f_w=0.
\]
Twisting again with parameter \(u\in\Cs\) multiplies \(we^z\) by
\(u\) and \(w^{-1}e^{-z}\) by \(u^{-1}\). This gives
\[
 (f_w)_u=1+(uw)e^z=f_{uw},\qquad f_{uw}''+(uw)^{-1}e^{-z}f_{uw}'-f_{uw}=0,
\]
the solution and equation obtained by twisting once with parameter \(uw\).
For \(r>1\), the same composition uses the coordinatewise product
\(\mathbf u\mathbf w=(u_1w_1,\ldots,u_rw_r)\).

For every character \(\chi\), linear independence of the exponentials
gives \(f_\chi\not\equiv0\). If \(z_0\) were a multiple zero of
\(f_\chi\), then \(f_\chi(z_0)=f_\chi'(z_0)=0\). Since \(A_\chi\)
and \(B_\chi\) are entire, uniqueness for the initial value problem in
\eqref{eq:twist} would give \(f_\chi\equiv0\). Thus every zero of
\(f_\chi\) is simple. We use this property when we pass to a limit
that retains only the terms on an edge of the Newton polytope.

\subsection{The convex hull in integer coordinates}
Let \(\mathcal{Q}_f=\{Q_1,\ldots,Q_m\}\) be the set of exponents of
\(f\) in \eqref{eq:sum}. If \(\#\mathcal{Q}_f=1\), then \(f=Pe^Q\)
and the degree bound is immediate. In this case the zeros of \(f\)
are the finitely many zeros of \(P\). We thus assume \(\#\mathcal{Q}_f\ge2\).

For \(\gamma\in\mathcal{Q}_f\), write \(P_\gamma\) for the multiplier
of \(e^\gamma\) in \(f\).
In the chosen \(\Z\)-basis \(\gamma_1,\ldots,\gamma_r\) of
\(\Gamma\), write the exponents and their coordinate vectors as
\[
 Q_j=\sum_{i=1}^r m_{ji}\gamma_i,
 \qquad p_j=(m_{j1},\ldots,m_{jr})\in\Z^r.
\]
The convex hull \(N=\conv\{p_1,\ldots,p_m\}\)
is a \emph{lattice polytope} \cite[Lecture~II, Sect.~1]{cox}.
Let \(X_1,\ldots,X_r\) be independent formal variables, with \(X_i\)
corresponding to \(e^{\gamma_i}\). Then \(N\) is the \emph{Newton polytope}
\cite[Definition~3.4]{maclagan} of the Laurent polynomial
\[
 \sum_{j=1}^m P_j(z)X_1^{m_{j1}}\cdots X_r^{m_{jr}}.
\]
For any \(\gamma\in\Gamma\), denote its integer coordinate vector in
this basis by \(p_\gamma\), so \(p_{Q_j}=p_j\). We use these coordinates
when viewing elements of \(\Gamma\) as lattice points.
A face of \(N\) is the set on which a real linear functional attains
its minimum. An edge is a one-dimensional face.

\noindent\textbf{Affine span.}\enspace
The \emph{linear span} \(\operatorname{span}_{\R}\{p_1,\ldots,p_m\}\)
is the set of all real linear combinations of these points. An affine subspace
is a translate of a linear subspace. The smallest affine subspace
containing the points is their \emph{affine span}, given by
\begin{align*}
 \operatorname{aff}\{p_1,\ldots,p_m\}
 &=\left\{\sum_{i=1}^m\lambda_i p_i :
    \lambda_i\in\R,\ \sum_{i=1}^m\lambda_i=1\right\}\\
 &=p_1+\operatorname{span}_{\R}\{p_2-p_1,\ldots,p_m-p_1\}.
\end{align*}
Its dimension is the dimension of the span of the \(m-1\) difference
vectors, so it is at most \(m-1\), with equality exactly when
\(p_2-p_1,\ldots,p_m-p_1\) are linearly independent.

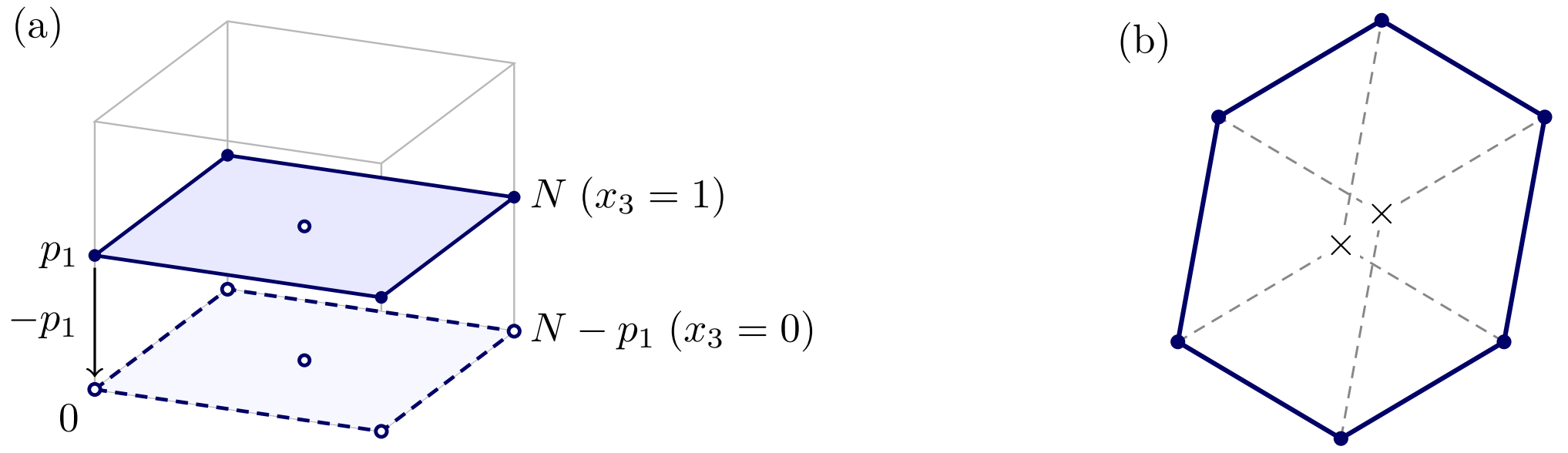
\begin{figure}[htbp]
  \centering
  \begin{minipage}[b]{0.49\linewidth}
    \centering
    \begin{tikzpicture}[
        x={(1.23cm,-0.18cm)}, y={(0.57cm,0.43cm)}, z={(0cm,1.15cm)},
        font=\small, hull/.style={draw=blue!40!black, line width=0.9pt},
        frame edge/.style={draw=black!28, line width=0.45pt}]
      % Ambient cube, and the two parallel positions of the lattice square.
      \draw[frame edge] (0,0,0) -- (2,0,0) -- (2,2,0) -- (0,2,0) -- cycle;
      \draw[frame edge] (0,0,2) -- (2,0,2) -- (2,2,2) -- (0,2,2) -- cycle;
      \foreach \vx/\vy in {0/0,2/0,2/2,0/2}{
        \draw[frame edge] (\vx,\vy,0) -- (\vx,\vy,2);
      }
      \path[hull, densely dashed, fill=blue!3]
        (0,0,0) -- (2,0,0) -- (2,2,0) -- (0,2,0) -- cycle;
      \path[hull, fill=blue!9]
        (0,0,1) -- (2,0,1) -- (2,2,1) -- (0,2,1) -- cycle;
      \draw[->, semithick, shorten >=3pt, shorten <=3pt]
        (0,0,1) -- node[midway,left] {$-p_1$} (0,0,0);
      \foreach \vx/\vy in {0/0,2/0,2/2,0/2}{
        \fill[blue!40!black] (\vx,\vy,1) circle[radius=1.55pt];
        \draw[hull, fill=white] (\vx,\vy,0) circle[radius=1.35pt];
      }
      \draw[hull, fill=white] (1,1,1) circle[radius=1.25pt];
      \draw[hull, fill=white] (1,1,0) circle[radius=1.25pt];
      \node[left] at (0,0,1) {$p_1$};
      \node[below left] at (0,0,0) {$0$};
      \node[right] at (2,2,1) {$N\;(x_3=1)$};
      \node[right] at (2,2,0) {$N-p_1\;(x_3=0)$};
      \node[anchor=north west] at (-0.85cm,3.42cm) {(a)};
    \end{tikzpicture}
  \end{minipage}\hfill
  \begin{minipage}[b]{0.49\linewidth}
    \centering
    \begin{tikzpicture}[
        x={(1.40cm,-0.83cm)}, y={(-1.40cm,-0.83cm)}, z={(0.35cm,1.93cm)},
        font=\small, remaining/.style={draw=blue!40!black, line width=1.15pt},
        removed/.style={draw=black!47, densely dashed, line width=0.55pt}]
      % View near a body diagonal so the six surviving edges form a visible cycle.
      \draw[removed] (1,0,0) -- (0,0,0) -- (0,1,0)
        (0,0,0) -- (0,0,1);
      \draw[removed] (1,1,0) -- (1,1,1) -- (0,1,1)
        (1,1,1) -- (1,0,1);
      \draw[remaining]
        (1,0,0) -- (1,1,0) -- (0,1,0) -- (0,1,1) --
        (0,0,1) -- (1,0,1) -- cycle;
      \foreach \vx/\vy/\vz in {1/0/0,1/1/0,0/1/0,0/1/1,0/0/1,1/0/1}{
        \fill[blue!40!black] (\vx,\vy,\vz) circle[radius=1.8pt];
      }
      \foreach \vx/\vy/\vz in {0/0/0,1/1/1}{
        \node[fill=white, inner sep=0.75pt] at (\vx,\vy,\vz) {$\times$};
      }
      \node[anchor=north west] at (-2.06cm,2.05cm) {(b)};
    \end{tikzpicture}
  \end{minipage}
  \caption{Affine coordinates and edge connectivity (AI-assisted).}
  \label{fig:affine-connectivity}
\end{figure}

The convex hull \(N\) has the same affine span as \(p_1,\ldots,p_m\).
If that span has dimension \(d\), translating by \(-p_1\) and choosing
a basis of its direction space identifies \(N\) with a full-dimensional
polytope in \(\R^d\). This affine bijection preserves vertices and edges.
Figure~\ref{fig:affine-connectivity}(a) shows a lattice square \(N\)
translated by \(-p_1\) from \(x_3=1\) to \(x_3=0\), where its first
two coordinates give a full-dimensional square in \(\R^2\).

Balinski's theorem states that the edge graph of a \(d\)-dimensional
polytope, \(d\ge1\), remains connected after the removal of fewer than
\(d\) vertices and their incident edges \cite[p.~432]{balinski}.
In Figure~\ref{fig:affine-connectivity}(b), deleting two opposite
vertices and their dashed edges leaves a cycle through the six remaining
vertices. We use only the consequence that any two vertices of \(N\)
can be joined by an edge path.

\FloatBarrier

We now argue by contradiction. Assume that
\[
 \deg(\gamma-\delta)\ge2
 \qquad\text{for some }\gamma,\delta\in\mathcal{Q}_f.
\]

For some \(j\ge2\), the coefficients of \(z^j\) in \(\gamma\) and
\(\delta\) differ. For \(x=(x_1,\ldots,x_r)\in\R^r\), let \(L_j(x)\)
be the coefficient of \(z^j\) in \(\sum_{k=1}^r x_k\gamma_k\).
The map \(L_j\colon\R^r\to\C\) is real-linear, so for a convex
combination of the vertices \(v_i\) of \(N\) we have
\[
 x=\sum_i\lambda_i v_i,\quad \lambda_i\ge0,\quad \sum_i\lambda_i=1
 \qquad\Longrightarrow\qquad L_j(x)=\sum_i\lambda_i L_j(v_i).
\]
Equal values at all vertices would make \(L_j\) constant on \(N\).
The distinct values at the points corresponding to \(\gamma\) and \(\delta\)
therefore force distinct values at two vertices.
By Balinski's theorem, an edge path joins these two vertices.
Since their \(L_j\)-values differ, some edge \(E\) along the path has
endpoints \(p_{\gamma_0},p_{\gamma_*}\) with different \(L_j\)-values.
Thus \(\deg(\gamma_*-\gamma_0)\ge j\ge2\).

Let \(l\ge1\) be the
greatest common divisor of the integer coordinates of \(\gamma_*-\gamma_0\)
and set \(u=(\gamma_*-\gamma_0)/l\). The coordinates of \(u\) have greatest
common divisor one, so \(u\) is called \emph{primitive}. As a polynomial,
\(u\) has the same degree as \(\gamma_*-\gamma_0\).

The points of \(E\) are \(p_{\gamma_0}+sp_u\) for \(0\le s\le l\).
By B\'ezout's identity \cite[Proposition~11, p.~91]{coppel}, the
coordinates of \(p_u\) have an integer linear combination equal to one.
Its coefficients define an integer-valued linear functional \(b\) with
\(b(p_u)=1\). We write \(b(\gamma)=b(p_\gamma)\).
For a lattice point \(p=p_{\gamma_0}+sp_u\), we have
\(s=b(p-p_{\gamma_0})\in\Z\). Thus the lattice points on \(E\) are
\(p_{\gamma_0}+jp_u\) for integers \(0\le j\le l\).
For the edge from \((1,0)\) to \((5,6)\), \(l=\gcd(4,6)=2\) and
\(p_u=(4,6)/2=(2,3)\), with lattice points \((1,0),(3,3),(5,6)\).
Here \(b(m_1,m_2)=-m_1+m_2\) gives \(1\) at \(p_u\), and \(2b\) gives
\(2\). We use a unit increase for simplicity.

\noindent\textbf{Supporting functional.}\enspace
A \emph{supporting functional} for a face of \(N\) is a real-linear map
\(a\colon\R^r\to\R\) whose minimum on \(N\) is attained exactly on that
face. We apply this definition to the edge \(E\).
For \(\gamma\in\Gamma\), we write \(a(\gamma)=a(p_\gamma)\), so the
minimum is \(k=a(\gamma_0)\). Writing \(a(x)=\sum_{i=1}^r c_i x_i\)
for \(x=(x_1,\ldots,x_r)\), we seek integer coefficients \(c_i\) so that
\(a(\Z^r)\subseteq\Z\).

If \(E=N\), take \(a=0\). Otherwise, the definition of a face gives
a real-linear functional \(a\) whose set of minimizers on \(N\) is
exactly \(E\). Thus \(a(x)=k=a(\gamma_0)\) for every \(x\in E\), and
\(a(\gamma)>k\) for every vertex \(p_\gamma\) of \(N\) outside \(E\).
Since \(\gamma_*-\gamma_0=lu\), the equal endpoint values give
\(l\,a(u)=a(\gamma_*)-a(\gamma_0)=0\). Hence the coefficient vector
\(c=(c_1,\ldots,c_r)\) is perpendicular to \(u\), and \(a(x)=c\cdot x\)
measures position in this normal direction.
The edge direction \(u=(u_1,\ldots,u_r)\) is fixed. Given a real normal
vector \(c\) and \(\varepsilon>0\), we seek \(\widetilde c\in\mathbb Q^r\)
with \(\widetilde c\cdot u=0\) and \(\|\widetilde c-c\|_2<\varepsilon\).
Approximate all but one coordinate of \(c\) by rationals
\(\widetilde c_j\), leaving a coordinate \(c_i\) unknown with \(u_i\ne0\).
To impose \(\widetilde c\cdot u=0\), set
$$\widetilde c_i=-u_i^{-1}\sum_{j\ne i}u_j\widetilde c_j.$$
This last coordinate is rational and tends to \(c_i\) as the other
coordinates \(\widetilde c_j\) tend to \(c_j\). Thus the rational normal
vectors are dense among the real vectors perpendicular to \(u\).
For the supporting normal \(c\), take
\(\varepsilon\) small enough that \(\widetilde c\cdot(p_\gamma-p_{\gamma_0})>0\)
for every vertex \(p_\gamma\) of \(N\) not on \(E\). Multiply \(\widetilde c\)
by a positive common denominator to obtain an integer normal vector.

\noindent\textbf{Normalized twist.}\enspace
Let \(t,w\in\Cs\) be two complex parameters. The map
\(\gamma\mapsto t^{a(\gamma)}w^{b(\gamma)}\) is a character of \(\Gamma\)
because \(a\) and \(b\) are additive. Applying the twist formula
\eqref{eq:twist-definition} to \(f\) with this character and dividing
the result by \(t^k\) gives the \emph{normalized twist}
\begin{equation}\label{eq:family}
 h_t(z,w)=\sum_{\gamma\in\mathcal{Q}_f}
 t^{a(\gamma)-k}w^{b(\gamma)}P_\gamma(z)e^{\gamma(z)}.
\end{equation}
For each \(\gamma\in\mathcal{Q}_f\) there are two cases.
If \(p_\gamma\in E\), then \(a(\gamma)-k=0\) and the factor
\(t^{a(\gamma)-k}\) is identically \(1\). If \(p_\gamma\notin E\),
then \(a(\gamma)-k\) is a positive integer and \(t^{a(\gamma)-k}\to0\)
as \(t\to0\). Thus the sum extends holomorphically to \(t=0\) for
\(w\ne0\), and only the terms on \(E\) remain at \(t=0\).
The convergence is uniform on compact subsets of \(\C\times\Cs\)
in \((z,w)\). For \(t,w\in\Cs\), every zero in \(z\) is simple because
\(h_t\) is a nonzero constant multiple of a twisted solution.

In the Laurent polynomial defining \(N\), the term with exponent vector
\(p_\gamma\) has \emph{weight} \(a(p_\gamma)\). The sum of the terms of
minimum weight is the \emph{initial form} with respect to \(a\)
\cite[Sect.~5]{maclagan}. The minimum here is \(k\), so both the initial
form and the limit \(h_0\) retain exactly the terms with \(p_\gamma\in E\).

For the sum \(1+e^z+e^{z^2}\) in the basis \(z,z^2\), the exponent
vectors are \((0,0),(1,0),(0,1)\), and the Laurent polynomial is
\(1+X_1+X_2\). To retain a quadratic exponent difference, choose the
vertical edge from \((0,0)\) to \((0,1)\), whose endpoint exponents
differ by \(z^2\). Take \(a(m_1,m_2)=m_1\), with minimum \(k=0\) on
this edge, and \(b(m_1,m_2)=m_2\). The corresponding torus element is
\((t,w)\). In each row of Table~\ref{tab:edge-limits}, the chosen edge
is the set of points of the triangle where \(a(m_1,m_2)=k\).

\begin{table}[htbp]
  \caption{Initial forms and limits for the edges of the triangle.}
  \label{tab:edge-limits}
  \centering
  \renewcommand{\arraystretch}{1.25}
  \setlength{\tabcolsep}{5pt}
  \begin{tabular}{@{}lcccccc@{}}
    \toprule
    Edge & \(a(m_1,m_2)\) & \(k=\min_N a\) & \(b(m_1,m_2)\) &
    Initial form & \(h_t(z,w)\) & \(h_0(z,w)\) \\
    \midrule
    Vertical & \(m_1\) & \(0\) & \(m_2\) & \(1+X_2\) &
    \(1+te^z+we^{z^2}\) & \(1+we^{z^2}\) \\
    Horizontal & \(m_2\) & \(0\) & \(m_1\) & \(1+X_1\) &
    \(1+we^z+te^{z^2}\) & \(1+we^z\) \\
    Diagonal & \(-m_1-m_2\) & \(-1\) & \(m_2\) & \(X_1+X_2\) &
    \(t+e^z+we^{z^2}\) & \(e^z+we^{z^2}\) \\
    \bottomrule
  \end{tabular}
\end{table}

\begin{samepage}
The vertical, horizontal, and diagonal edges have endpoint exponent
differences \(z^2\), \(z\), and \(z^2-z\), respectively. Thus the
vertical and diagonal edges both give a difference of degree at least
two. The edge labels refer to the chosen lattice basis. Changing that
basis changes the coordinates used to describe the exponent polynomials.
Expressing the same functionals \(a,b\) in the new coordinates preserves
their values on each exponent and hence gives the same \(h_t\) and \(h_0\).
\par
\end{samepage}

\begin{figure}[!b]
  \centering
  \begin{tikzpicture}[
    x={(0.75cm,0.55cm)}, y={(1.80cm,0cm)}, z={(0cm,1.35cm)},
    font=\small,
    cube/.style={draw=black!50, line width=0.5pt},
    hidden/.style={cube, densely dashed},
    minimum/.style={draw=blue!40!black, line width=1.2pt}]
    \foreach \offset/\panel/\element/\constraint in
      {0/a/{(t,t,w)}/{m_1=m_2=0},5.2/b/{(t,1,w)}/{m_1=0}} {
      \begin{scope}[xshift=\offset cm]
        \draw[hidden] (0,0,0) -- (1,0,0) -- (1,1,0)
          (1,0,0) -- (1,0,1);
        \if\panel b
          \fill[blue!12, fill opacity=0.8]
            (0,0,0) -- (0,1,0) -- (0,1,1) -- (0,0,1) -- cycle;
        \fi
        \draw[cube] (0,0,0) -- (0,1,0) -- (0,1,1) -- (0,0,1) -- cycle
          (0,1,0) -- (1,1,0) -- (1,1,1) -- (1,0,1) -- (0,0,1)
          (0,1,1) -- (1,1,1);
        \foreach \i/\j/\k in
          {0/0/0,0/0/1,0/1/0,0/1/1,1/0/0,1/0/1,1/1/0,1/1/1}{
            \fill[black!48] (\i,\j,\k) circle[radius=1.3pt];
          }
        \if\panel a
          \draw[minimum] (0,0,0) -- (0,0,1);
          \foreach \k in {0,1} {
            \fill[blue!40!black] (0,0,\k) circle[radius=1.9pt];
          }
        \else
          \draw[minimum] (0,0,0) -- (0,1,0) -- (0,1,1) -- (0,0,1) -- cycle;
          \foreach \j/\k in {0/0,0/1,1/0,1/1}{
            \fill[blue!40!black] (0,\j,\k) circle[radius=1.9pt];
          }
        \fi
        \node at (1.28cm,2.20cm) {(\panel) \(\element\)};
        \node at (1.28cm,-0.32cm) {\(\constraint\)};
      \end{scope}
    }
    \begin{scope}[shift={(3.50cm,0.49cm)}, inner sep=1.5pt]
      \draw[->] (0,0) -- (0.42cm,0.31cm) node[above right] {\(m_1\)};
      \draw[->] (0,0) -- (0.66cm,0) node[right] {\(m_2\)};
      \draw[->] (0,0) -- (0,0.52cm) node[above] {\(m_3\)};
    \end{scope}
  \end{tikzpicture}
  \caption{Minimum sets of the two \(t\)-weights (AI-assisted).}
  \label{fig:cube-minima}
\end{figure}
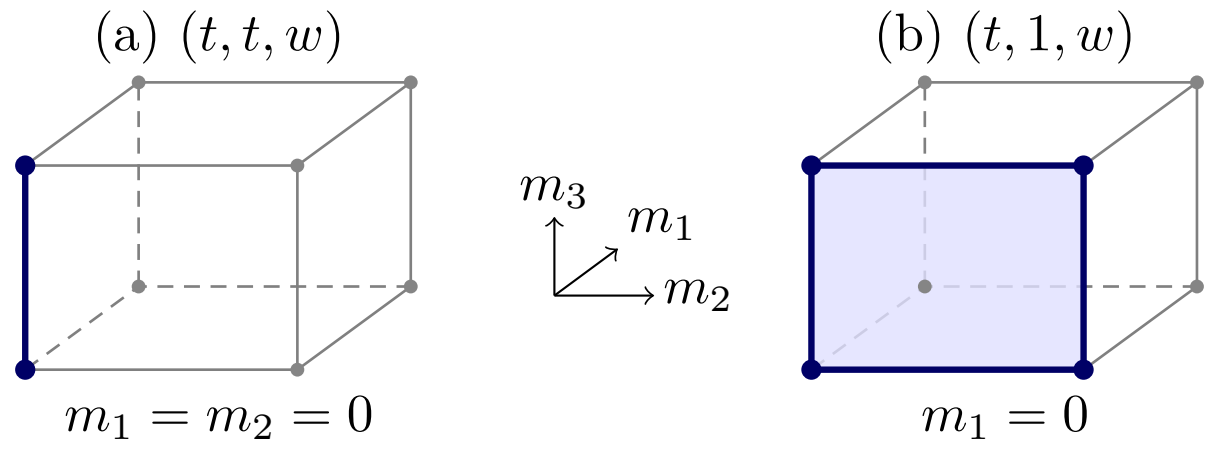

\begin{samepage}
Suppose that the exponent vectors are the eight vertices of the cube
\([0,1]^3\). The torus element \((t,t,w)\) varies its first two
coordinates with the same \(t\). It multiplies the term at
\((m_1,m_2,m_3)\) by \(t^{m_1+m_2}w^{m_3}\). Since \(m_1,m_2\ge0\),
the power of \(t\) is minimized exactly where \(m_1=m_2=0\). The
minimum set is the edge \(\{0\}\times\{0\}\times[0,1]\) in
Figure~\ref{fig:cube-minima}(a). As \(t\to0\), its two vertex terms
remain with factors \(1\) and \(w\).
\par
\end{samepage}

\FloatBarrier

\begin{samepage}
For the same cube, \((t,1,w)\) gives the factor \(t^{m_1}w^{m_3}\).
The supporting functional \(m_1\) is minimized exactly on the square face
\(\{0\}\times[0,1]^2\) in Figure~\ref{fig:cube-minima}(b). Its four
vertex terms remain, two with factor \(1\) and two with \(w\).
\par
\end{samepage}

At \(t=0\) with \(w\ne0\) fixed, \eqref{eq:family} retains the terms on
\(E\). Let \(p_u\) be the lattice direction vector of the polynomial
\(u\). A step along \(E\) gives
\[
\begin{array}{@{}ll@{}}
 \text{In polynomial notation} & \gamma\longmapsto\gamma+u,\\[2pt]
 \text{In lattice coordinates} & p_\gamma\longmapsto p_\gamma+p_u.
\end{array}
\]
The endpoint difference \(lp_u\) represents \(l\) lattice steps.
Its coordinates have gcd \(l\), so those of \(p_u\) are coprime and
B\'ezout's identity gives \(b(p_u)=1\).
Thus \(b(p_\gamma+p_u)=b(p_\gamma)+1\), and each step multiplies
\(w^{b(\gamma)}e^{\gamma(z)}\) by \(we^{u(z)}\), giving \((we^{u(z)})^j\)
after \(j\) steps. Setting \(Q=u\) gives
\begin{equation}\label{eq:edge}
 h_0(z,w)=w^{b(\gamma_0)}e^{\gamma_0(z)}P(z,we^{Q(z)}),
 \qquad
 P(z,x)=\sum_{\substack{0\le j\le l\\\gamma_0+ju\in\mathcal{Q}_f}}P_{\gamma_0+ju}(z)x^j .
\end{equation}
Both endpoint exponents occur in \(f\), so the \(x^0\) and \(x^l\) coefficients are nonzero polynomials.

\section{Branching in the limiting equation}\label{sec:branching}

For fixed $w\ne0$, the zeros of $h_0(\cdot,w)$ satisfy $P(z,we^{Q(z)})=0$,
since the prefactor in \eqref{eq:edge} is nonvanishing. For the vertical edge
in Table~\ref{tab:edge-limits}, $h_0=1+we^{z^2}$. Setting $x=we^{z^2}$
gives $P(z,x)=1+x$, so we take $H=P$ and obtain
\[
 H(z,x)=0\quad\Longrightarrow\quad x=-1\quad\Longrightarrow\quad
 w=xe^{-z^2}=-e^{-z^2}.
\]
Near $z=0$, write $w+1=z^2v(z)$ with $v(0)=1$ and set $s=z\sqrt{v(z)}$
with $\sqrt{v(0)}=1$, so $s^2=w+1$. For small $\varepsilon>0$, as $\theta$
increases from $0$ to $2\pi$, $w=-1+\varepsilon e^{i\theta}$ circles $-1$
once, while the roots $s=\pm\sqrt{\varepsilon}e^{i\theta/2}$ each rotate
through $\pi$ about $0$. Since $s'(0)=1$, the local inverse carries
the exchange to the two zeros $z_\pm$ near $0$, where $z_-=-z_+$.

\subsection{The surface of an algebraic function}\label{subsec:algebraic-function}

To study the zeros in \eqref{eq:edge}, we apply the following theorem
to a monic irreducible factor $H$ of $P(z,x)$ over $\C(z)$.

\smallskip
\noindent\textbf{Algebraic function theorem}
\textup{\cite[Theorem~8.9]{forster}.}
Let $H\in\C(z)[x]$ be monic and irreducible with $d=\deg_x H\ge1$.
Its coefficients $c_j(z)$ are rational in $z$ and meromorphic on $\widehat{\C}$.
There exist a connected Riemann surface $\Sigma$, a proper
$d$-sheeted holomorphic branched covering
$z\colon\Sigma\to\widehat{\C}$ of the Riemann sphere, and a
meromorphic function $x$ on $\Sigma$ such that
\[
 H\bigl(z(p),x(p)\bigr)
 = x(p)^d+c_1\bigl(z(p)\bigr)x(p)^{d-1}+\cdots+c_d\bigl(z(p)\bigr)=0
\]
as an identity of meromorphic functions on $\Sigma$.
 Since $z$ is proper and $\widehat{\C}$ is compact,
$\Sigma=z^{-1}(\widehat{\C})$ is compact.

Fix $z_0$ at which no coefficient of $H$ has a pole. If $H(z_0,x)=0$
has $d$ distinct roots in $x$, the map $p\mapsto x(p)$ gives a
bijection from $z^{-1}(z_0)$ to these roots. For $H=1+x$, there is
one sheet and $x\equiv-1$ on $\Sigma$.

For $H=x^2-z$, take $\Sigma=\widehat{\C}$ with affine coordinate $x$
and projection $z=x^2$. The inverse image of $U=\{|z-1|<1/2\}$
consists of two neighborhoods on this same $\Sigma$, one about $x=-1$
and one about $x=1$. Each maps biholomorphically onto $U$
(Figure~\ref{fig:square-cover}).

\begingroup
\setlength{\intextsep}{6pt plus 2pt minus 1pt}
\begin{figure}[!ht]
  \centering
  \begin{tikzpicture}[
      font=\small,
      sheet/.style={fill=blue!8, draw=blue!40!black, line width=0.8pt},
      disk/.style={fill=black!4, draw=black!45, line width=0.65pt},
      projection/.style={->, draw=blue!40!black, line width=0.8pt},
      declare function={
        rootreal(\a)=sqrt((sqrt(1.25+cos(\a))+1+0.5*cos(\a))/2);}]
    \path[draw=black!38, fill=black!2, line width=0.5pt]
      (-2.42,1.18)
      .. controls (-2.42,1.51) and (-1.45,1.62) .. (-0.58,1.59)
      .. controls (0.18,1.56) and (1.06,1.68) .. (1.80,1.54)
      .. controls (2.24,1.46) and (2.42,1.35) .. (2.42,1.17)
      .. controls (2.42,.86) and (1.49,.76) .. (.52,.79)
      .. controls (-.20,.82) and (-1.00,.69) .. (-1.74,.83)
      .. controls (-2.14,.91) and (-2.42,1.02) .. cycle;
    \node[anchor=east] at (-2.53,1.21) {\(\Sigma\)};
    \node[anchor=east] at (-2.53,.63) {\(z=x^2\)};
    \node[anchor=east] at (-2.53,0) {\(\widehat{\C}\)};
    \begin{scope}[shift={(0,1.21)}, x=1.15cm, y=1.15cm]
      % The two square roots of 1 + (1/2)e^{i a} trace the boundaries.
      \foreach \sgn in {-1,1} {
        \path[sheet] plot[variable=\ang,domain=0:360,samples=81]
          ({\sgn*rootreal(\ang)},{\sgn*0.25*sin(\ang)/rootreal(\ang)})
          -- cycle;
        \fill[blue!40!black] (\sgn,0) circle[radius=1.45pt];
      }
    \end{scope}
    \node[above, inner sep=1pt] at (-1.15,1.64) {\(U_-\;(x=-1)\)};
    \node[above, inner sep=1pt] at (1.15,1.64) {\(U_+\;(x=1)\)};
    \draw[projection] (-1.10,.92) to[out=-66,in=153] (-.33,.33);
    \draw[projection] (1.10,.92) to[out=-114,in=27] (.33,.33);
    \path[disk] (0,0) circle[radius=.43];
    \fill[blue!40!black] (0,0) circle[radius=1.45pt];
    \node[anchor=east,font=\footnotesize] at (-.55,0) {\(z=1\)};
    \node[anchor=west] at (.64,0) {\(U=\{|z-1|<\tfrac12\}\)};
  \end{tikzpicture}
  \caption{Two sheets of \(z=x^2\) over a disk about \(1\) (AI-assisted).}
  \label{fig:square-cover}
\end{figure}
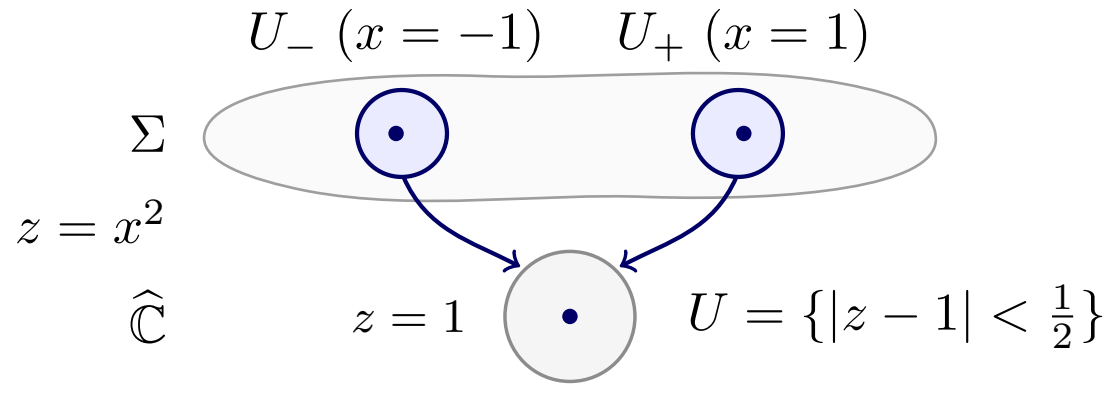
\FloatBarrier
\endgroup

For $z_0\in\C$ at which all coefficients of $H$ are finite, let
$x_1(z_0),\ldots,x_d(z_0)$ be the roots of $H(z_0,x)=0$, listed with
multiplicity. Set
\[
 \Delta_H(z_0)=\prod_{i<j}\bigl(x_i(z_0)-x_j(z_0)\bigr)^2 .
\]
This product is $1$ when $d=1$ and vanishes exactly when a root is repeated.
It is symmetric in the roots and hence is a polynomial in the
coefficients $c_j$. Thus $\Delta_H$ is rational in $z$ and is not
identically zero by the irreducibility of $H$ over $\C(z)$.
For a monic quadratic in $x$, Vieta's formulas give
$\Delta_H(z)=c_1(z)^2-4c_2(z)$. In the example $H=x^2-z$ of
Figure~\ref{fig:square-cover}, this gives
$
 \Delta_H(z)=4z.
$

Define the finite subset $D_H$ of the $z$-sphere by
\begin{equation}\label{eq:excluded-z-values}
 D_H=\{\infty\}
 \cup\bigcup_{j=1}^{d}\bigl\{z_0\in\C:c_j\text{ has a pole at }z_0\bigr\}
 \cup\bigl\{z_0\in\C:\Delta_H(z_0)=0\bigr\}.
\end{equation}
For each $z_0\in D_H$, there are at most $d$ points of $\Sigma$ with
$z=z_0$, so $z^{-1}(D_H)$ is finite. At any $p\in\Sigma$, choose a
local coordinate $s$ with $s(p)=0$ and shrink its disk to exclude all
points of $z^{-1}(D_H)$ other than possibly $p$. On the punctured disk
$0<|s|<\varepsilon$, distinct points with the same $z$ represent
distinct $x$-roots. Thus $s\mapsto(z(s),x(s))$ is one-to-one there.

\subsection{Counting zeros and poles of a differential}

A \emph{meromorphic differential} on $\Sigma$ is written locally as
$F(s)\,ds$, with $F$ meromorphic. Its order at a point is the
order of $F$ in a coordinate $s$ centered there, positive at a zero
and negative at a pole. This order is independent of the coordinate,
since a change of local coordinate has nonvanishing derivative.
For a nonzero meromorphic differential $\Lambda$ on a compact surface
of genus $g$, we use \cite[Theorem~17.12]{forster}
\begin{equation}\label{eq:differential-order-sum}
 \sum_{p\in\Sigma}\operatorname{ord}_p\Lambda=2g-2.
\end{equation}

\subsection{A branch point of the limiting equation}

By \eqref{eq:edge}, the limiting equation $h_0(z,w)=0$ is equivalent to
$P(z,we^{Q(z)})=0$ for $w\ne0$.
As in $w+1=s^2$, the next lemma gives a local parametrization with
several distinct $z$ values for each $w\ne w_0$ sufficiently close
to the branch value $w_0$.

\begin{lemma}\label{lem:branch}
Suppose that $P\in\C[z,x]$ has at least two distinct powers of $x$
with nonzero coefficients, and $Q\in\C[z]$ has degree at least two.
Then there are holomorphic functions $z(s)$ and $w(s)$ near $s=0$ with
$w(0)\ne0$ and $P(z(s),w(s)e^{Q(z(s))})=0$. On a sufficiently small
punctured disk, $s\mapsto(z(s),w(s))$ is a one-to-one parametrization
of its image. Namely, for $s_1,s_2$ in that disk,
\[
 (z(s_1),w(s_1))=(z(s_2),w(s_2))
 \quad\Longrightarrow\quad s_1=s_2.
\]
Moreover, $w(s)-w(0)$ has a zero of finite order at least two at $s=0$.
\end{lemma}

This order is called the \emph{local degree} of $s\mapsto w(s)$ at $0$.

\begin{proof}
Regard $P$ as a polynomial in $x$ over $\C(z)$. After removing its
lowest power of $x$, we get a polynomial of positive degree with
nonzero constant coefficient. Choose a monic irreducible factor
$H$ of this polynomial and apply the algebraic function theorem of
Section~\ref{subsec:algebraic-function}. On the resulting compact
surface $\Sigma$, the relation $P(z,x)=0$ holds as an identity of
meromorphic functions, and $x$ is not identically zero.

At points with $z$ finite and $x$ finite and nonzero, we want to
find a critical point of $w=xe^{-Q(z)}$. Its logarithmic derivative
extends to all of $\Sigma$ as the meromorphic differential
\begin{equation}\label{eq:eta}
 \eta=\frac{dx}{x}-d(Q(z)).
\end{equation}
Use a local coordinate $s$ centered at the point under consideration
for both terms in \eqref{eq:eta}. Writing $x(s)=s^m u(s)$ with $u$
holomorphic and $u(0)\ne0$ gives $dx/x=(m/s+u'/u)\,ds$. The poles of
$dx/x$ are simple and occur precisely at the zeros and poles of $x$.

At a pole of $z$ of order $k$, we have $z(s)=a s^{-k}+\cdots$ with
$a\ne0$. The leading term of $Q(z(s))$ is a nonzero multiple of
$s^{-k\deg Q}$, so $d(Q(z))$ has a pole of order $k\deg Q+1\ge3$.
This pole persists in $\eta$ since $dx/x$ has at most a simple pole.
The covering map $z$ is a nonconstant meromorphic function on the
compact surface $\Sigma$. If $z$ had no poles, it would be holomorphic
on all of $\Sigma$ and $|z|$ would attain its maximum there. The maximum
modulus principle would force $z$ to be constant. Thus $z$ has at least
one pole, and $\eta$ has a pole of order at least three.

By \eqref{eq:differential-order-sum}, the sum of the orders of
$\eta$ is $2g-2\ge-2$, so a pole of order at least three forces at
least one zero $p$. At a zero or pole
of $x$ where $z$ is finite, $d(Q(z))$ is holomorphic, so the simple
pole of $dx/x$ persists in $\eta$. Thus $z(p)\in\C$ and $x(p)\in\Cs$,
so $w=xe^{-Q(z)}$ is holomorphic and nonvanishing near $p$.

Take a local coordinate $s$ centered at $p$ and write $\eta=F(s)\,ds$.
Since $\eta$ vanishes at $p$, we have $F(0)=0$, and $dw=w\eta$ gives
$w'(0)=w(0)F(0)=0$. Thus the Taylor expansion is
\[
 w(s)-w(0)=\frac{w''(0)}{2!}s^2+\frac{w'''(0)}{3!}s^3+\cdots .
\]
If all its coefficients were zero, $w$ would be constant near $p$
and $\eta=dw/w$ would vanish on that neighborhood. The identity
theorem would then give $\eta\equiv0$ on $\Sigma$, contradicting
the existence of a pole of $\eta$ over $z=\infty$. Hence there is
a least integer $k\ge2$ for which $w^{(k)}(0)\ne0$, and $w-w(p)$
vanishes to finite order $k$ at $p$.

On the punctured disk of Section~\ref{subsec:algebraic-function},
$z(s)$ avoids $D_H$ in \eqref{eq:excluded-z-values}, so
$s\mapsto(z(s),x(s))$ is one-to-one.
Two parameters with the same $(z,w)$ have the same $(z,x)$ by
$x=we^{Q(z)}$, so they coincide. Substitution also gives
$P(z(s),w(s)e^{Q(z(s))})=P(z(s),x(s))=0$.
\end{proof}

\subsection{Completion of the proof}

\begin{proof}[Proof of Theorem~\ref{thm:exponents}]
Assume two exponents differ by a polynomial of degree at least two
and take $h_t$ from \eqref{eq:family}. Lemma~\ref{lem:branch} gives a
branch $(z(s),w(s))$ of $h_0(z,w)=0$ with $h_0$ as in \eqref{eq:edge}.
The point $p\in\Sigma$ chosen in the proof of Lemma~\ref{lem:branch}
corresponds to $s=0$, where $w'(0)=0$, and gives the branch point
$(z_0,w_0)=\bigl(z(p),x(p)e^{-Q(z(p))}\bigr)$.
We will show that every zero near $(z_0,w_0)$ has $z=g_j(w)$ for some
$1\le j\le M$, where $g_1,\ldots,g_M$ are holomorphic near $w_0$.
In \eqref{eq:edge}, the exponents $\gamma_0+ju$ are distinct and
$w_0^{b(\gamma_0)}P_{\gamma_0}\not\equiv0$ since $w_0\ne0$. By linear
independence, $h_0(\cdot,w_0)\not\equiv0$. Since its zero $z_0$ is
isolated, choose a disk $U$ centered at $z_0$ so small that $z_0$
is the only zero of $h_0(\cdot,w_0)$ in $\overline U$.
By continuity, choose a disk $V$ about $w_0$ with $\overline V\subset\Cs$
and $h_0\ne0$ on $\partial U\times\overline V$.

By \eqref{eq:family}, $h_t\to h_0$ uniformly on
$\overline U\times\overline V$ as $t\to0$. Choose $\varepsilon>0$
so that $h_t\ne0$ on $\partial U\times\overline V$ for
$|t|<\varepsilon$. For each such $t$ and each $w\in V$, the function
$h_t(\cdot,w)$ is entire by \eqref{eq:family}. The argument principle
gives the number of its zeros in $U$, counted with multiplicity, as
\[
 M(t,w)=\frac{1}{2\pi i}\int_{\partial U}
          \frac{\partial_z h_t(z,w)}{h_t(z,w)}\,dz.
\]
By \eqref{eq:family}, the numerator and denominator depend continuously
on $(z,t,w)$ on the contour, where the denominator is nonzero. Thus
$M(t,w)$ is continuous and integer-valued on the connected set
$\{|t|<\varepsilon\}\times V$, so it has a constant value $M$.
Taking $s=0$ in Lemma~\ref{lem:branch} gives $h_0(z_0,w_0)=0$
by \eqref{eq:edge}. By the choice of $U$, $M=M(0,w_0)$ is the
multiplicity of this zero, so $M\ge1$.

Fix a sufficiently small $t\ne0$. For each $w\in V$, the normalized
twist $h_t(\cdot,w)$ in \eqref{eq:family} is a nonzero constant multiple
of the corresponding nonzero twisted solution $f_\chi$, so it solves
the second-order linear equation \eqref{eq:twist}. Its zeros are simple
by uniqueness for the initial value problem. Thus the count
$M(t,w)=M$ gives $M$ distinct zeros in $U$ for every $w\in V$.

Let $Z_t=\{(z,w)\in U\times V:h_t(z,w)=0\}$ and let
$\pi\colon Z_t\to V$ be the projection onto the $w$-coordinate,
$\pi(z,w)=w$. For any $w_*\in V$, denote these distinct zeros
of $h_t(\cdot,w_*)$ by $z_1,\ldots,z_M$. Choose pairwise disjoint disks
$D_j\subset U$ about $z_j$. Since $\partial_z h_t(z_j,w_*)\ne0$,
the holomorphic implicit function theorem gives, after shrinking
the disks if needed, a common disk $V_*\subset V$ about $w_*$ such that the
zeros in each $D_j\times V_*$ are given by $z=z_{t,j}(w)$.
Here $z_{t,j}$ is holomorphic on $V_*$ and $z_{t,j}(w_*)=z_j$.
On each of these graphs, $w\mapsto(z_{t,j}(w),w)$ is the holomorphic
local inverse of $\pi$, so $\pi$ is a local biholomorphism.

For any compact $K\subset V$, the absence of zeros on
$\partial U\times K$ gives
\[
 \pi^{-1}(K)
 =\{(z,w)\in\overline U\times K:h_t(z,w)=0\}.
\]
This is a closed subset of the compact set $\overline U\times K$
by continuity of $h_t$, so $\pi$ is proper.

For each fixed $w\in V_*$, the $M$ local graphs give $M$ distinct
simple zeros, one in each $D_j$. Since $M(t,w)=M$, these are all
the zeros of $h_t(\cdot,w)$ in $U$. Thus the inverse image is the
disjoint union of the local graphs
\[
 \begin{aligned}
 Z_t\cap(D_j\times V_*)&=\{(z_{t,j}(w),w):w\in V_*\},\\
 \pi^{-1}(V_*)&=\bigsqcup_{j=1}^{M}\bigl(Z_t\cap(D_j\times V_*)\bigr).
 \end{aligned}
\]
Here $\bigsqcup$ denotes disjoint union.
The projection maps each graph biholomorphically onto $V_*$.
Since every $w_*\in V$ has such a neighborhood $V_*$, the projection
$\pi\colon Z_t\to V$ is an $M$-sheeted covering.
As $V$ is simply connected, the root functions extend holomorphically
to $z_{t,1},\ldots,z_{t,M}\colon V\to U$
\cite[Theorems~4.17 and~4.22]{forster}. Thus for each sufficiently
small $t\ne0$ the whole zero set over $V$ is
\[
 Z_t=\bigsqcup_{j=1}^{M}\{(z_{t,j}(w),w):w\in V\}.
\]

\enlargethispage{1.5\baselineskip}
\noindent
\begin{minipage}[t]{.645\linewidth}
\vspace{0pt}
\setlength{\parindent}{1em}
\noindent
Figure~\ref{fig:spiral-cover} uses $z=s(s-1)$ and $w=1+sz$.
The critical values $c_0=1$ and $c_1=23/27$ lie in the hole of $A$.
At $w=c_0$, $s^2(s-1)=0$. Over $V$, the double root $s=0$ gives two
blue sheets ($m=2$), and the simple root $s=1$ gives one orange sheet
($M=3$ in total). The cut is $L$. Pale gold shading marks $z(s)\in U$.

Take $t_n\to0$ with $t_n\ne0$. Apply Montel's theorem to the uniformly
bounded branches $z_{t_n,j}$ in turn, taking a further subsequence
each time. The final subsequence preserves all the locally uniform
limits $g_j$, which are holomorphic on $V$. Fix $w\in V$. If the $g_j(w)$ miss
a zero of $h_0(\cdot,w)$, take a closed disk $K\subset U$ about it
that avoids the $g_j(w)$ and has no zeros on its boundary. For small
$t$, uniform convergence gives $|h_t-h_0|<|h_0|$ on that boundary.
By Rouch\'e's theorem the zero counts in $K$, with multiplicity,
are equal, so $h_t$ has a zero there. Yet along the subsequence
$z_{t_n,j}(w)\to g_j(w)\notin K$ for each $j$, so eventually $K$
contains no roots.
\end{minipage}\hfill
\begin{minipage}[t]{.32\linewidth}
  \vspace{0pt}
  \centering
  \includegraphics[width=.96\linewidth,height=.31\textheight,keepaspectratio]{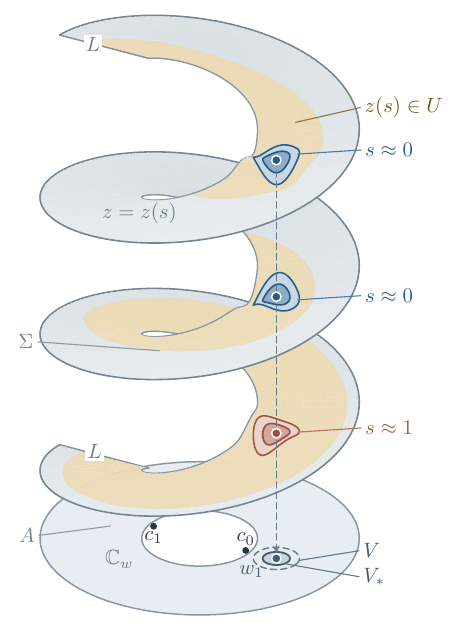}\par
  \nointerlineskip
  \setlength{\captionindent}{0pt}
  \setlength{\abovecaptionskip}{2pt}
  \captionof{figure}{Sheets (AI aid).}
  \label{fig:spiral-cover}
\end{minipage}
\par\addvspace{10pt}\noindent
Thus $h_0(z,w)=0$ in $U\times V$ implies $z=g_j(w)$ for some $j$.
For each small $s$, applying this to $(z(s),w(s))$ gives an index $j$
that may depend on $s$, so
\[
\prod_{j=1}^M\bigl(z(s)-g_j(w(s))\bigr)=0.
\]
The factors are holomorphic in $s$, so the identity theorem gives a
fixed $j_0$ such that
\begin{equation}\label{eq:fixed-limit-root}
 z(s)=g_{j_0}(w(s))\qquad\text{for all sufficiently small }s.
\end{equation}
By Lemma~\ref{lem:branch}, the local degree of $w(s)$ at $0$ is $m\ge2$.
For nearby $w\ne w_0$, the local normal form
\cite[Theorem~2.1 and Remark~2.2]{forster} gives $m$ distinct nonzero
preimages $s_1,\ldots,s_m$ (see Figure~\ref{fig:square-cover} for the square map).
By the injectivity in Lemma~\ref{lem:branch}, the $z(s_i)$ are distinct zeros
of $h_0(\cdot,w)$ in $U$, so $2\le m\le M$. Yet $w(s_1)=w(s_2)=w$,
and \eqref{eq:fixed-limit-root} gives $z(s_1)=g_{j_0}(w)=z(s_2)$,
a contradiction.

The zero bound follows from Dickson's theorem as explained after
the statement of Theorem~\ref{thm:exponents}.
\end{proof}

\begin{remark}
For fixed small $t\ne0$, $z_{t,j}(w_0)\mapsto z_{t,j}(w)$ gives a
holomorphic motion of the zeros in $U$ over $V$ \cite[p.~185]{slodkowski}.
Along the chosen subsequence $t_n\to0$,
$z_{t_n,i}-z_{t_n,j}\to g_i-g_j$ locally uniformly on $V$.
For $i\ne j$ the differences on the left have no zeros, so Hurwitz's
theorem gives $g_i\equiv g_j$ or $g_i(w)\ne g_j(w)$ throughout $V$.
After merging identical limits, the zeros of $h_0(\cdot,w)$ in $U$
would also move holomorphically. At $w_0$ the set is $\{z_0\}$, so its
motion gives one point per $w$, whereas Lemma~\ref{lem:branch} gives
at least two distinct zeros for nearby $w\ne w_0$, a contradiction.
\end{remark}

\section{Proof of Corollary~\ref{thm:hitw}}\label{sec:hitw}
\begin{proof}
By \cite[Theorem~4.16(a)]{survey}, $q=\rho(f)=\rho(A)$ is a positive
integer, $B$ is an exponential polynomial of order at most $q-1$, and
\[
 f(z)=c+\sum_{j=1}^m F_j(z)e^{w_jz^q},
 \qquad c\ne0.
\]
The $w_j$ are distinct and nonzero, and the nonzero $F_j$ are
exponential polynomials of order at most $q-1$. Since $B$ is an
exponential polynomial, Theorem~\ref{thm:exponents} applies.
On expanding $F_j$,
each exponent in its summand has the form $w_jz^q+R(z)$ with
$\deg R<q$. Its difference from $0$, the exponent of the nonzero
constant term $c$, has degree $q$. Theorem~\ref{thm:exponents} now gives
$q=1$.
An exponential polynomial of order at most zero is a polynomial.
\end{proof}

\begin{samepage}
The polynomial multipliers in \eqref{eq:sum} matter here. As observed
in \cite[(3.14)]{dual}, the function $1+e^{z^2}$ satisfies
\eqref{eq:ode} with
\[
 A(z)=\frac{e^{-z^2}-1}{2z}-2z,\qquad B(z)=-1.
\]
Although the singularity at zero is removable, the multiplier of
$e^{-z^2}$ is $1/(2z)$. Uniqueness of the representation over
$\C(z)$ precludes a representation of $A$ with polynomial multipliers.
\end{samepage}

\section{Equations of prescribed order}\label{sec:prescribed}

Every exponential polynomial satisfies a monic linear differential equation
with polynomial coefficients \cite{vpt}. For monic equations of prescribed order \(n\)
with entire coefficients, Gol'dberg's criterion is that all zeros of
\(f\not\equiv0\) have multiplicity less than \(n\) \cite[Sect.~4.1]{survey}.
Twisting preserves monic equations with exponential polynomial coefficients,
so the following criterion tests every twist \(f_\chi\) from
\eqref{eq:twist-definition}, where \(\chi\) ranges over the characters of the
group \(\Gamma_f\) generated by the normalized exponents of \(f\).

\begin{proposition}\label{prop:twist-order}
Let \(f\) be a nonzero exponential polynomial and \(n\ge1\) an integer.
There exist exponential polynomials \(A_0,\ldots,A_{n-1}\) such that
\[
 f^{(n)}+\sum_{j=0}^{n-1}A_j f^{(j)}=0
\]
if and only if every zero of every twist \(f_\chi\), for \(\chi\) a
character of \(\Gamma_f\), has multiplicity less than \(n\).
\end{proposition}

\begin{proof}
Suppose first that such an equation is given. For each \(j\), let
\(\widetilde A_j\) be the sum of the terms of \(A_j\) whose normalized
exponents lie in \(\Gamma_f\), the group of integer combinations of
the normalized exponents of \(f\). Differentiation changes only the
polynomial multipliers,
so all exponents of \(f^{(j)}\) lie in \(\Gamma_f\). When a term of
\(A_j\) with exponent \(\delta\) is multiplied by a term of \(f^{(j)}\)
with exponent \(\gamma\), the resulting exponent \(\delta+\gamma\)
lies in \(\Gamma_f\) exactly when \(\delta\) does, since
\(\gamma\in\Gamma_f\) and \(\delta=(\delta+\gamma)-\gamma\).
Thus the terms of \(A_j f^{(j)}\) with exponents in \(\Gamma_f\) sum
to \(\widetilde A_j f^{(j)}\).

By linear independence of distinct normalized exponentials, each
polynomial multiplier in the original identity vanishes after grouping
terms with the same exponent. Its terms with exponents in \(\Gamma_f\)
therefore give
\[
 f^{(n)}+\sum_{j=0}^{n-1}\widetilde A_j f^{(j)}=0.
\]

For example, take \(f=e^z\), so \(\Gamma_f=\Z z\). Substitution gives
\[
 f''+(e^{z^2}-1)f'-e^{z^2}f=(1-1)e^z+(1-1)e^{z^2+z}=0.
\]
The two groups have exponents \(z\in\Gamma_f\) and \(z^2+z\notin\Gamma_f\).
Discarding the second group removes the \(e^{z^2}\)-terms in the
coefficients and leaves \(f''-f'=0\).

In general, twisting preserves products and commutes with
differentiation. Thus for any character \(\chi\) of \(\Gamma_f\), the
equation with coefficients \(\widetilde A_j\) becomes
\[
 (f_\chi)^{(n)}+
 \sum_{j=0}^{n-1}(\widetilde A_j)_\chi (f_\chi)^{(j)}=0.
\]
The leading coefficient stays \(1\) since \(\chi(0)=1\), and the other
coefficients are exponential polynomials, hence entire. Since
\(f_\chi\not\equiv0\), a zero of multiplicity at least \(n\) would
give zero initial data, contradicting uniqueness.

Conversely, choose a basis \(\eta_1,\ldots,\eta_s\) of \(\Gamma_f\).
Let \(R=\C[z,X_1^{\pm1},\ldots,X_s^{\pm1}]\) be the Laurent polynomial
ring in \(X_1,\ldots,X_s\) over \(\C[z]\). For \(0\le j<n\), obtain
\(H_j(z,X)\in R\) from \(f^{(j)}\) by replacing each \(e^\gamma\) by
\(\prod_i X_i^{m_i}\), where \(\gamma=\sum_i m_i\eta_i\).
For example, when \(n\ge2\), take the basis \(\eta_1=z\), \(\eta_2=z^2\).
The exponent \(z^2-z=\eta_2-\eta_1\) gives \(X_1^{-1}X_2\).
Thus
\[
 \begin{aligned}
 f=1+ze^{z^2-z}+e^z
 &\quad\longleftrightarrow\quad H_0=1+zX_1^{-1}X_2+X_1,\\
 f'=(1+2z^2-z)e^{z^2-z}+e^z
 &\quad\longleftrightarrow\quad H_1=(1+2z^2-z)X_1^{-1}X_2+X_1.
 \end{aligned}
\]
For any \(a\in\C\) and nonzero \(x_1,\ldots,x_s\), choose the character with
\(\chi(\eta_i)=x_i e^{-\eta_i(a)}\). Iterating the differentiation rule in
Section~\ref{subsec:twisting} gives \((f_\chi)^{(j)}=(f^{(j)})_\chi\).
For a term \(P(z)e^{\gamma(z)}\) of \(f^{(j)}\), write
\(\gamma=\sum_i m_i\eta_i\). Both \(\chi\) and \(\gamma\mapsto e^{\gamma(a)}\)
send sums to products. By \eqref{eq:twist-definition}, the twisted term at \(a\) is
\[
 P(a)\chi(\gamma)e^{\gamma(a)}
 =P(a)\prod_i\bigl(\chi(\eta_i)e^{\eta_i(a)}\bigr)^{m_i}
 =P(a)\prod_i x_i^{m_i}.
\]
This is the corresponding term in \(H_j(a,x_1,\ldots,x_s)\). Summing gives
\[
 H_j(a,x_1,\ldots,x_s)=(f_\chi)^{(j)}(a),\qquad 0\le j<n.
\]
For the example above, the values on the two sides are
\[
 \begin{aligned}
 H_0(a,x_1,x_2)&=1+ax_1^{-1}x_2+x_1=f_\chi(a),\\
 H_1(a,x_1,x_2)&=(1+2a^2-a)x_1^{-1}x_2+x_1=(f_\chi)'(a).
 \end{aligned}
\]
Suppose all of \(H_0,\ldots,H_{n-1}\) vanish at the same point
\((a,x_1,\ldots,x_s)\in\C\times(\Cs)^s\). The single character
\(\chi\) chosen at this point would give \((f_\chi)^{(j)}(a)=0\) for
every \(0\le j<n\), so \(f_\chi\) would have a zero of multiplicity
at least \(n\) at \(a\). This contradicts the hypothesis, so the
\(H_j\) have no common zero in \(\C\times(\Cs)^s\).

Let \(\widehat H_j\) be the polynomial obtained from \(H_j\) by replacing
\(X_i^{-k}\) with \(Y_i^k\) for \(k>0\). A common zero of
\(\widehat H_0,\ldots,\widehat H_{n-1}\) and of \(X_iY_i-1\) for
\(i=1,\ldots,s\) would satisfy \(X_i\ne0\) and \(Y_i=X_i^{-1}\) for every \(i\).
It would therefore give a common zero of \(H_j\) for \(j=0,\ldots,n-1\).
Since these polynomials have no common zero, the weak Nullstellensatz
\cite[Ch.~4, \S1, Theorem~1]{clo} gives polynomials \(C_j,D_i\) in
\(z,X_1,\ldots,X_s,Y_1,\ldots,Y_s\) such that
\[
 \sum_{j=0}^{n-1} C_j\widehat H_j + \sum_{i=1}^s D_i(X_iY_i-1)=1.
\]
Substituting \(Y_i=X_i^{-1}\) sends \(\widehat H_j\) to \(H_j\) and
\(C_j\) to \(U_j\). Both \(H_j\) and \(U_j\) belong to \(R\).
Each factor \(X_iY_i-1\) becomes \(X_iX_i^{-1}-1=0\), so
\[
 \sum_{j=0}^{n-1}U_jH_j=1\qquad\text{in }R.
\]
Now substitute \(X_i=e^{\eta_i(z)}\) and set
\(u_j(z)=U_j(z,e^{\eta_1(z)},\ldots,e^{\eta_s(z)})\).
This substitution corresponds to taking \(\chi(\eta_i)=1\) for every
\(i\), so \(H_j\) evaluates to the original \(f^{(j)}\).
Thus \(\sum_{j=0}^{n-1}u_jf^{(j)}=1\), with each \(u_j\) an exponential
polynomial.
Since \(f^{(n)}\) is again an exponential polynomial, so is each product
\(A_j=-f^{(n)}u_j\). Multiplying the identity by \(-f^{(n)}\) gives
\[
 f^{(n)}+\sum_{j=0}^{n-1} A_j f^{(j)}=0.\qedhere
\]
\end{proof}

\begin{samepage}
For example, let \(f=1+e^Q\) with \(Q\) a nonconstant normalized
polynomial. The twists are \(1+we^Q\) for \(w\in\Cs\). At a point
\(a\), only the twist with \(w=-e^{-Q(a)}\) vanishes. Let \(m=\deg_a Q\)
be the local degree of \(Q\) at \(a\). For some \(c\ne0\), the Taylor expansions give
\[
 \begin{aligned}
 Q(z)-Q(a)&\sim c(z-a)^m,\\
 Q'(z)&\sim mc(z-a)^{m-1},\\
 1+we^{Q(z)}=1-e^{Q(z)-Q(a)}&\sim -c(z-a)^m
 \end{aligned}
 \qquad (z\to a).
\]
The zero has multiplicity \(m=1+\operatorname{ord}_a Q'\), where
\(\operatorname{ord}_aQ'=0\) if \(Q'(a)\ne0\). By Proposition~\ref{prop:twist-order},
the least order of a monic equation with exponential polynomial coefficients is
\[
 1+\max_{a\in\C}\deg_a Q
 =2+\max_{a\in\C}\operatorname{ord}_a Q'.
\]
\end{samepage}
This gives the first two rows of Table~\ref{tab:prescribed-orders},
since the corresponding derivatives \(Q'\) are \(qz^{q-1}\) and
\(z^{q-1}-1\), respectively. For the last row, the twists are
\((1-we^z)^k\), \(w\in\Cs\), and each zero has multiplicity \(k\).

\par\medskip
\noindent\begin{minipage}{\linewidth}
  \captionof{table}{Exponent degrees and least equation orders.}
  \label{tab:prescribed-orders}
  \centering
  \smallskip
  \renewcommand{\arraystretch}{1.2}
  \begin{tabular}{@{}lccc@{}}
    \toprule
    \(f\) & \shortstack{Largest exponent\\degree}
      & \shortstack{Maximal zero multiplicity\\among twists}
      & \shortstack{Least equation\\order} \\
    \midrule
    \(1+e^{z^q}\), \(q\ge1\) & \(q\) & \(q\) & \(q+1\) \\
    \(1+e^{z^q/q-z}\), \(q\ge2\) & \(q\) & \(2\) & \(3\) \\
    \((1-e^z)^k\), \(k\ge1\) & \(1\) & \(k\) & \(k+1\) \\
    \bottomrule
  \end{tabular}
\end{minipage}
\par\medskip

Using row~2 of Table~\ref{tab:prescribed-orders}, we construct a dominant coefficient of \(f'\) dual to \(f\).

\begin{proof}[Proof of Proposition~\ref{prop:higher}]
Let \(Q(z)=z^q/q-z\), \(P=Q'\), and \(f=1+e^Q\). For \(j\ge1\),
put \(R_j=e^{-Q}f^{(j)}\). These are polynomials, since \(R_1=P\)
and \(R_{j+1}=R_j'+PR_j\). In particular, \(R_2=P'+P^2\).
At any root \(a\) of \(R_1=P=z^{q-1}-1\), we have \(a\ne0\) and
\[
 R_2(a)=P'(a)+P(a)^2=P'(a)=(q-1)a^{q-2}\ne0.
\]
Thus \(R_1,R_2\) are coprime, so B\'ezout gives
\(U_1R_2+V_1R_1=1\) with \(U_1,V_1\in\C[z]\).
For \(n\ge3\), \(R_n=UR_2+VR_1\) (with \(U=R_nU_1\), \(V=R_nV_1\)), hence
\(f^{(n)}-Uf''-Vf'=0\).
Adding \((1+e^{-Q})f'-Pf=0\) gives
\[
 f^{(n)}-Uf''+(1+e^{-Q}-V)f'-Pf=0.
\]
The functions \(f\) and \(A_1=1+e^{-Q}-V\) both have growth order
\(q\), with nonzero highest-order frequencies \(1/q\) and \(-1/q\),
respectively. Every other coefficient is a polynomial.
\end{proof}

\clearpage
\section*{Declaration of AI assistance}
The author used \href{https://openai.com/codex/}{OpenAI Codex} to search
the literature, develop and check proofs, draft and revise the text and
\LaTeX\ source, prepare figures, and develop and test the Lean code.
The author takes responsibility for the mathematical content, references,
figures, and presentation.

\end{document}